# Constrained similarity of 2-D trajectories by minimizing the $H^1$ semi-norm of the trajectory difference

**Stefan M. Filipov [a], Atanas V. Atanasov [a], Ivan D. Gospodinov [a*]**

[a] Department of Computer Science,
University of Chemical Technology and Metallurgy,
8 "Kl. Ohridski" Blvd.,
Sofia, 1756, Bulgaria

* Corresponding author
email: idg2@cornell.edu

**Abstract**

This paper defines constrained functional similarity between 2-D trajectories via minimizing the $H^1$ semi-norm of the difference between the trajectories. An exact general solution is obtained for the case wherein the components of the trajectories are mesh-functions defined on a uniform mesh and the imposed constraints are linear. Various examples are presented, one of which features application to mechanics and two-point boundary value problems. A MATLAB code is given for the solution of one of the examples. The code could easily be adjusted to other cases.

**Introduction**

Suppose a trajectory is given and a new trajectory is sought that meets a number of imposed constraints and is as similar in behaviour to the original trajectory as possible without necessarily being close [1] to it. Such shape optimisation problems may have wide range of applications in many engineering fields [2] such as mechanics, fluid mechanics, aerodynamics, general transport phenomena, design and engineering of machines and equipment, etc. In [3] the authors have introduced constrained functional similarity between real-valued functions of one real variable via minimizing the $H^1$ semi-norm [4] of the difference between the functions. An exact general solution for mesh-functions has been presented. The similarity of trajectories in two and more dimensions is as important. This work defines constrained similarity between 2-D trajectories and provides an exact solution to the discretized case. Application to mechanics and two-point boundary value problems is presented in Example 4 of the Results section.

**Constrained similarity of 2-D trajectories**

Let $r^*(t)=(x^*(t),y^*(t))$ and $r(t)=(x(t),y(t))$ be two radius vectors whose components are real-valued functions of a real independent variable $t\in[a,b]$. The functions $r^*$ and $r$ define two 2-D trajectories. The trajectory $r^*$ will be *similar* to $r$, under certain given constraints, if $r^*$ minimizes the square of the $H^1$ semi-norm of the difference $r^*-r$:

$$|r^*-r|_{H_1}^2=\int_a^b\left(\frac{dr^*}{dt}-\frac{dr}{dt}\right)^2dt=\int_a^b\left(\frac{dx^*}{dt}-\frac{dx}{dt}\right)^2dt+\int_a^b\left(\frac{dy^*}{dt}-\frac{dy}{dt}\right)^2dt \tag{1}$$

and at the same time satisfies the constrains in question. The constraints that $r^*$ satisfies must be linear in $x^*$ and $y^*$. For example, linear combinations of functional values $x^*(t_i)$ and $y^*(t_i)$ at certain points $t_i$, integral constraints like $\int_a^b f(t)x^*(t)dt=1$ or $\int_a^b g(t)y^*(t)dt=1$, etc.

**Exact solution for discretized trajectories under linear constraints**

Partitioning the interval $t\in[a,b]$ by $N$ mesh points into $N-1$ intervals of equal size defines a *uniform mesh* on the interval: $\{t_i=a+(i-1)h,\ i=1,2,\ldots,N,\ h=(b-a)/(N-1)\}$, where $h$ is the step-size of the mesh. Let the trajectory $r$ be defined on the mesh, i.e. $\{r_i=r(t_i),\ i=1,2,\ldots,N\}$. In order to define constrained similarity between the trajectories $r^*$ and $r$ expression (1) is discretized using the forward finite differences $(x^*_{i+1}-x^*_i)/h$, etc. for the respective derivatives $dx^*/dt$, etc. at $t_i$, $i=1,2,\ldots,N-1$ and the integral is replaced by a sum. The constant $h$ is omitted since constant factors do not affect the minimization. Thus, the following objective function is obtained:

$$I=\sum_{i=1}^{N-1}((x^*_{i+1}-x^*_i)-(x_{i+1}-x_i))^2+\sum_{i=1}^{N-1}((y^*_{i+1}-y^*_i)-(y_{i+1}-y_i))^2. \tag{2}$$

In order to use the formulas derived in [3] we denote $x_i=u_i$, $y_i=u_{N+i}$, $x^*_i=u^*_i$, and $y^*_i=u^*_{N+i}$, for $i=1,2,\ldots,N$ and introduce the two vectors $u=[x_1,\ldots,x_N,y_1,\ldots,y_N]^T$ and $u^*=[x^*_1,\ldots,x^*_N,y^*_1,\ldots,y^*_N]^T$. The minimum of $I$ is sought subject to $M$ linear constraints:

$$\sum_{i=1}^{2N}A_{ji}u^*_i=c_j,\qquad j=1,2,\ldots,M<2N\,. \tag{3}$$

The constraints (3) can be written in a matrix form as $Au^*=c$, where

$$A=\begin{bmatrix}A_{11} & A_{12} & . & . & . & A_{1(2N)}\\ A_{21} & A_{22} & . & . & . & A_{2(2N)}\\ . & . & . & . & . & .\\ A_{M1} & A_{M2} & . & . & . & A_{M(2N)}\end{bmatrix},\quad c=\begin{bmatrix}c_1\\ c_2\\ .\\ c_M\end{bmatrix} \tag{4}$$

and $u^*$ is the $2N\times1$ column-vector of the unknowns. To find the minimum of $I$ subject to constraints (3) the Lagrange's method of the undetermined coefficients [5] is used. First, the Lagrangian

$$J = I + \sum_{j=1}^{M}\left(\lambda_j (c_j - \sum_{i=1}^{2N} A_{ji} u^*_i)\right) \tag{5}$$

is introduced, where $\lambda_j$, $j=1,2,\ldots,M$ are the Lagrange's undetermined coefficients. Then, the derivatives of $J$ with respect to $u^*_k$, $k=1,2,\ldots,2N$ are equated to zero:

$$\begin{aligned}
&\frac{\partial J}{\partial u^*_k} = 2((u^*_k - u^*_{k-1}) - (u_k - u_{k-1})) - 2((u^*_{k+1} - u^*_k) - (u_{k+1} - u_k)) - \sum_{j=1}^{M} \lambda_j A_{jk} = 0, \quad k = 2,..,N-1, N+2,..,2N-1 \\
&\frac{\partial J}{\partial u^*_1} = -2((u^*_2 - u^*_1) - (u_2 - u_1)) - \sum_{j=1}^{M} \lambda_j A_{j1} = 0 \\
&\frac{\partial J}{\partial u^*_N} = 2((u^*_N - u^*_{N-1}) - (u_N - u_{N-1})) - \sum_{j=1}^{M} \lambda_j A_{jN} = 0 \\
&\frac{\partial J}{\partial u^*_{N+1}} = -2((u^*_{N+2} - u^*_{N+1}) - (u_{N+2} - u_{N+1})) - \sum_{j=1}^{M} \lambda_j A_{j(N+1)} = 0 \\
&\frac{\partial J}{\partial u^*_{2N}} = 2((u^*_{2N} - u^*_{2N-1}) - (u_{2N} - u_{2N-1})) - \sum_{j=1}^{M} \lambda_j A_{j(2N)} = 0
\end{aligned} \tag{6}$$

The system of equations (6) is rearranged so that only terms containing $u^*_i$ remain on the left-hand side. Then, the system is written in a matrix form as:

$$\bar{L}u^* = \bar{L}u - \frac{1}{2}A^T\lambda, \tag{7}$$

where $\lambda$ is the $M\times1$ column-vector of the undetermined coefficients and $\bar{L}$ is the $2N\times2N$ matrix:

$$\bar{L} = \begin{bmatrix}
-1 & 1 & 0 & 0 & \cdots & & & & & 0 \\
1 & -2 & 1 & 0 & \cdots & & & & & 0 \\
0 & 1 & -2 & 1 & \cdots & & & & & 0 \\
 & & \ddots & & & & & \ddots & & \\
0 & \cdots & 1 & -2 & 1 & \cdots & & & & 0 \\
0 & \cdots & 0 & 1 & -1 & \cdots & & & & 0 \\
0 & \cdots & & & & -1 & 1 & 0 & 0 \cdots & 0 \\
0 & \cdots & & & & 1 & -2 & 1 & 0 \cdots & 0 \\
0 & \cdots & & & & 0 & 1 & -2 & 1 \cdots & 0 \\
 & & \ddots & & & & & \ddots & & \\
0 & \cdots & & & & & & 1 & -2 & 1 \\
0 & \cdots & & & & & & 0 & 1 & -1
\end{bmatrix} \tag{8}$$

In order to remove the singularity of $\bar{L}$ at least one of the equations for the constraints needs to be added to one of the first $N$ equations in (7) and at least one of the equations for the constraints needs to be added to one of the second $N$ equations in (7). For this reason a $2N\times2N$ matrix $\bar{A}$ is introduced whose first row is any row of the matrix $A$ (say row $j$) that corresponds to an $x$-constraint and whose row $N+1$ is any row of the matrix $A$ (say row $m$) that corresponds to a $y$-constraint. The rest of the elements of $\bar{A}$ are zeros. A $2N\times1$ column-vector $\bar{c}$ with components zeros except for $\bar{c}(1) = c(j)$ and $\bar{c}(N+1) = c(m)$ is also introduced. If necessary, more equations from $A$ can be included in $\bar{A}$. Now, the results for $u^*$ and $\lambda$, derived in [3], can be used:

$$u^* = u - (\bar{L} + \bar{A})^{-1}\left(\frac{1}{2}A^T\lambda + \bar{A}u - \bar{c}\right), \tag{9}$$

$$\lambda = 2\left(A(\bar{L} + \bar{A})^{-1}A^T\right)^{-1}\left(Au - c - A(\bar{L} + \bar{A})^{-1}(\bar{A}u - \bar{c})\right), \tag{10}$$

where $\bar{L}$ is defined in (8) . The right-hand side of (10) contains only known quantities. Once the column vector $\lambda$ is calculated it is substituted into (9) and the sought $u^*_i$, $i=1,2,\ldots,2N$ are found.

## Results

In this paragraph several examples are presented with three types of constraints: boundary, difference, and integral constraints. The last example describes an application to mechanics and two-point boundary value problems.

### Example 1

Consider the trajectory $r$ defined by $\{r_i=(x_i,y_i),\ x_i=\sin(2t_i),\ y_i=\sin(3t_i),\ i=1,2,\ldots,N\}$ on a uniform mesh with boundaries $a=0$, $b=\pi$, and number of mesh-points $N=101$. Using (9) and (10), the trajectory $r^*$, i.e. $\{r^*_i=(x^*_i,y^*_i),\ i=1,2,\ldots,N\}$, similar to $r$ and satisfying the following two integral and one difference constraints

$$\sum_{i=1}^{N} x^*_i = \sum_{i=1}^{N} x_i + \Delta S_x,\quad \sum_{i=1}^{N} y^*_i = \sum_{i=1}^{N} y_i + \Delta S_y,\quad x^*_k - x^*_{N-k+1} = 0 \tag{11}$$

is found for several values of $\Delta S_x$, $\Delta S_y$, and $k$ (see fig.1 below).

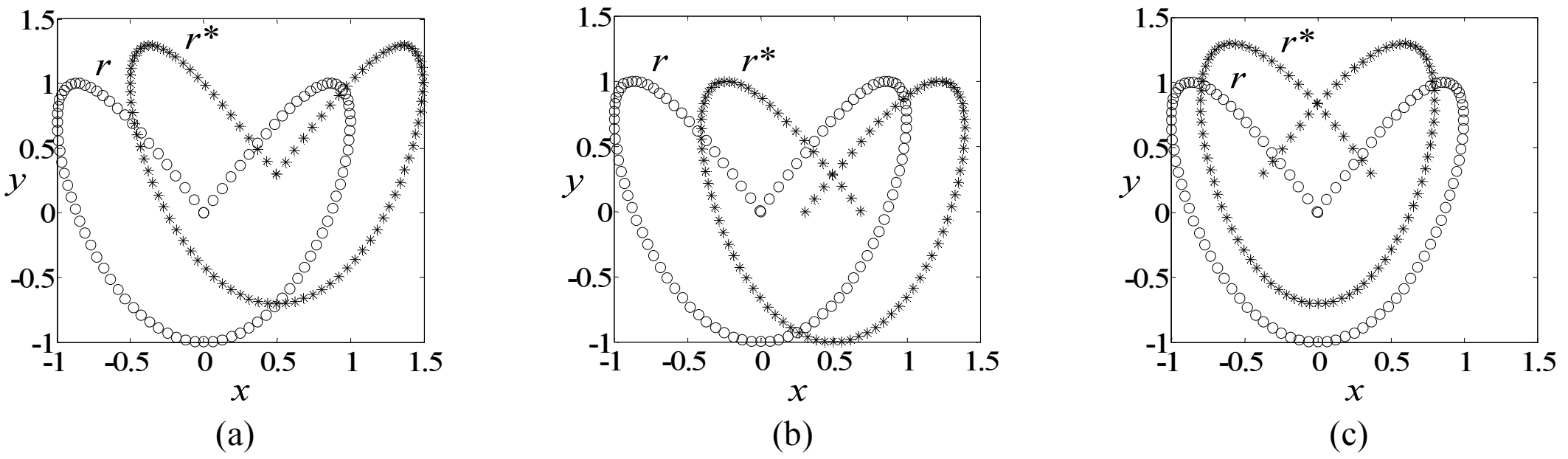

Fig.1. The original trajectory $r$ (o) and the similar to it trajectory $r^*$ (*) satisfying constraints (11) for (a) $\Delta S_x=50$, $\Delta S_y=30$, and $k=1$; (b) $\Delta S_x=50$, $\Delta S_y=0$, and $k=4$; and (c) $\Delta S_x=0$, $\Delta S_y=30$, and $k=7$.

### Example 2

Consider the trajectory $r$ defined by $\{r_i=(x_i,y_i),\ x_i= t_i-2\sin(t_i),\ y_i=1-2\cos(t_i),\ i=1,2,\ldots,N\}$ on a uniform mesh with $a=-\pi$, $b=3\pi$, and $N=101$. The trajectory $r^*$ similar to $r$ and satisfying the following boundary constraints

$$(x^*_1, y^*_1) = (x_1, y_1) + (\Delta x_1, \Delta y_1),\quad (x^*_N, y^*_N) = (x_N, y_N) + (\Delta x_N, \Delta y_N) \tag{12}$$

is found for several values of $(\Delta x_1, \Delta y_1)$ and $(\Delta x_N, \Delta y_N)$ (see fig 2. below).

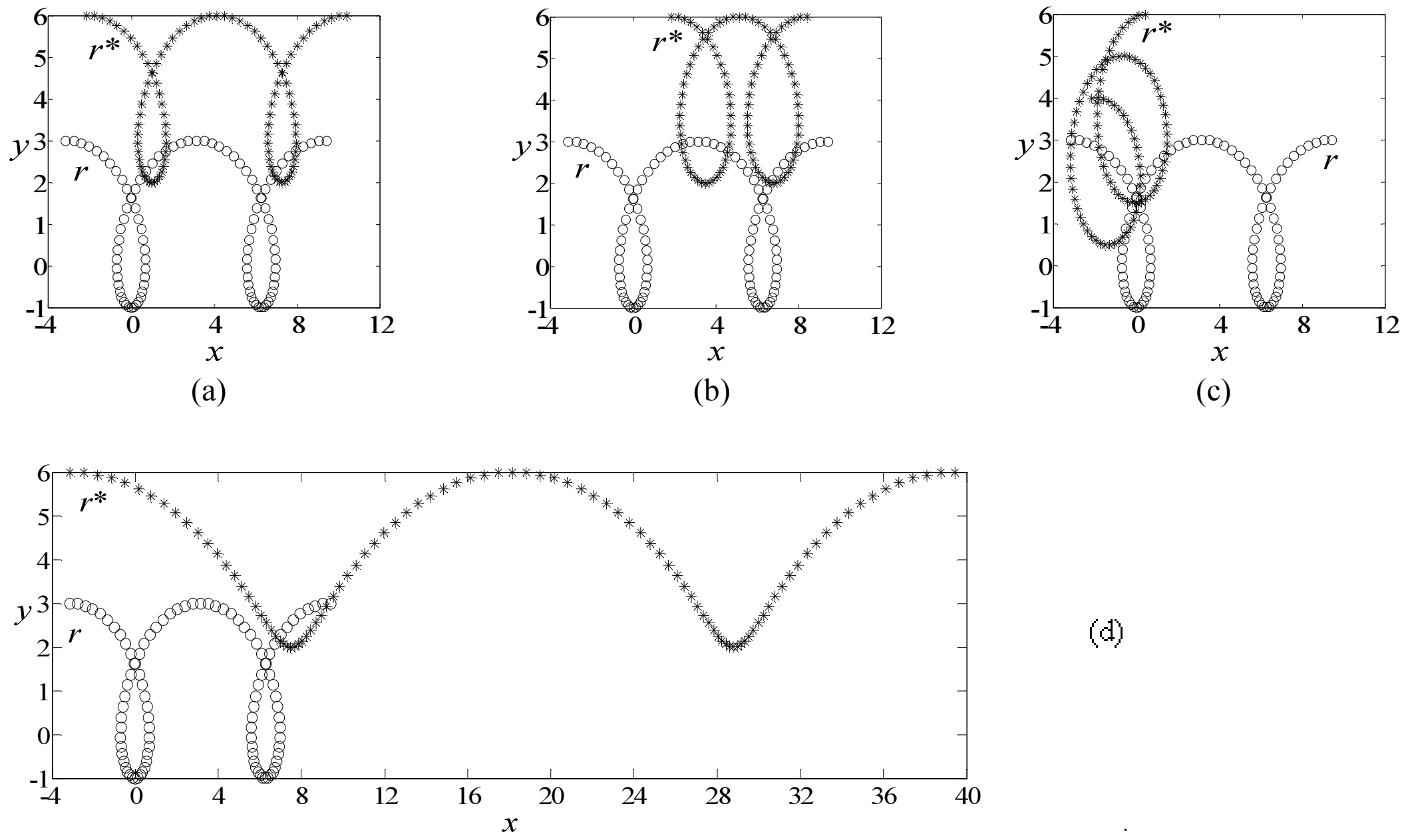

Fig. 2. The original trajectory $r$ (o) and the similar to it trajectory $r^*$ (*) satisfying constraints (12) for (a) $(\Delta x_1, \Delta y_1)=(1,3)$, $(\Delta x_N, \Delta y_N)=(1,3)$; (b) $(\Delta x_1, \Delta y_1)=(5,3)$, $(\Delta x_N, \Delta y_N)=(-1,3)$; (c) $(\Delta x_1, \Delta y_1)=(1,1)$, $(\Delta x_N, \Delta y_N)=(-9,3)$; and (d) $(\Delta x_1, \Delta y_1)=(0,3)$, $(\Delta x_N, \Delta y_N)=(30,3)$.

**Example 3**

Consider the trajectory $r$ defined by $\{r_i=(x_i,y_i),\ x_i=\sin(2t_i),\ y_i=(1-\sin(t_i))\sin(t_i),\ i=1,2,\ldots,N\}$ on a uniform mesh with $a=0$, $b=2\pi$, and $N=101$. The trajectory $r$ * similar to $r$ and satisfying the following four integral and two difference constraints:

$$\sum_{i=1}^{N} x^*_i = \sum_{i=1}^{N} x_i + 10,\ \sum_{i=1}^{N} y^*_i = \sum_{i=1}^{N} y_i,$$
$$\sum_{i=1}^{N} t_i x^*_i = \sum_{i=1}^{N} t_i x_i,\ \sum_{i=1}^{N} (t_i-a)(t_i-b) y^*_i = \sum_{i=1}^{N} (t_i-a)(t_i-b) y_i + \Delta T_y,$$
$$x^*_1 - x^*_N = 0,\ y^*_1 - y^*_N = 0 \tag{13}$$

is found for several values of $\Delta T_y$ (see fig.3 below).

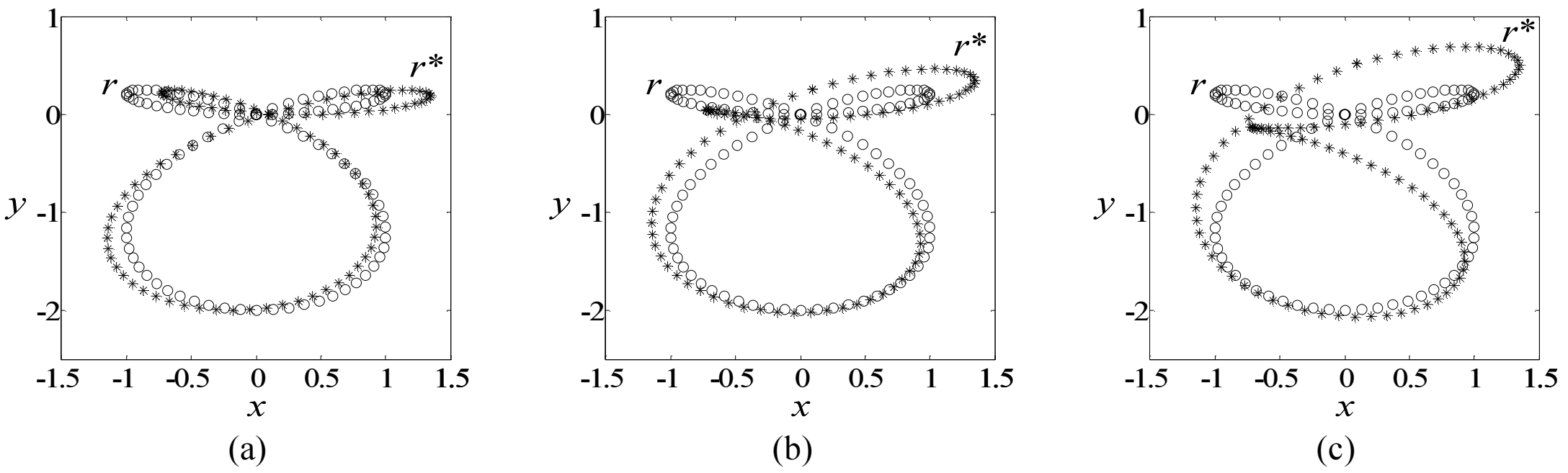

Fig.3. The original trajectory $r$ (o) and the similar to it trajectory $r^*$ (*) satisfying constraints (13) for (a) $\Delta T_y=0$; (b) $\Delta T_y=50$; and (c) $\Delta T_y=100$.

**Example 4**

A mass point, initially at rest, travels 2 seconds under the influence of the gravitational potential $U=gy$, $g=9.8$ (m/s$^2$). Placing the origin of the coordinate system at the initial position of the point and partitioning the time interval $t\in[0,2]$ by $N=11$ equally separated mesh-points the following discretized trajectory $r$ is obtained: $\{r_i=(x_i,y_i),\ x_i=0,\ y_i=-gt_i^2/2,\ t_i=(i-1)h,\ i=1,2,\ldots,N,\ h=0.2\ \text{(s)}\}$. The trajectory $r$ * similar to $r$ and satisfying the following boundary constraints:

$$(x^*_1, y^*_1) = (0,0),\quad (x^*_N, y^*_N) = (x_b, y_b), \tag{14}$$

is found for several values of $(x_b,y_b)$ (see fig.4 below). The obtained trajectory $r$ * describes *exactly* the motion of a point travelling for 2 seconds between points (0,0) and $(x_b,y_b)$ under the influence of the given potential. If the force field is not homogenous the trajectory $r$ * will describe the motion of the point only approximately. Then, however, $r$ * could be incorporated into a 'shooting-projection' iterative procedure to obtain the exact solution to the two-point boundary value problem [6], [7].

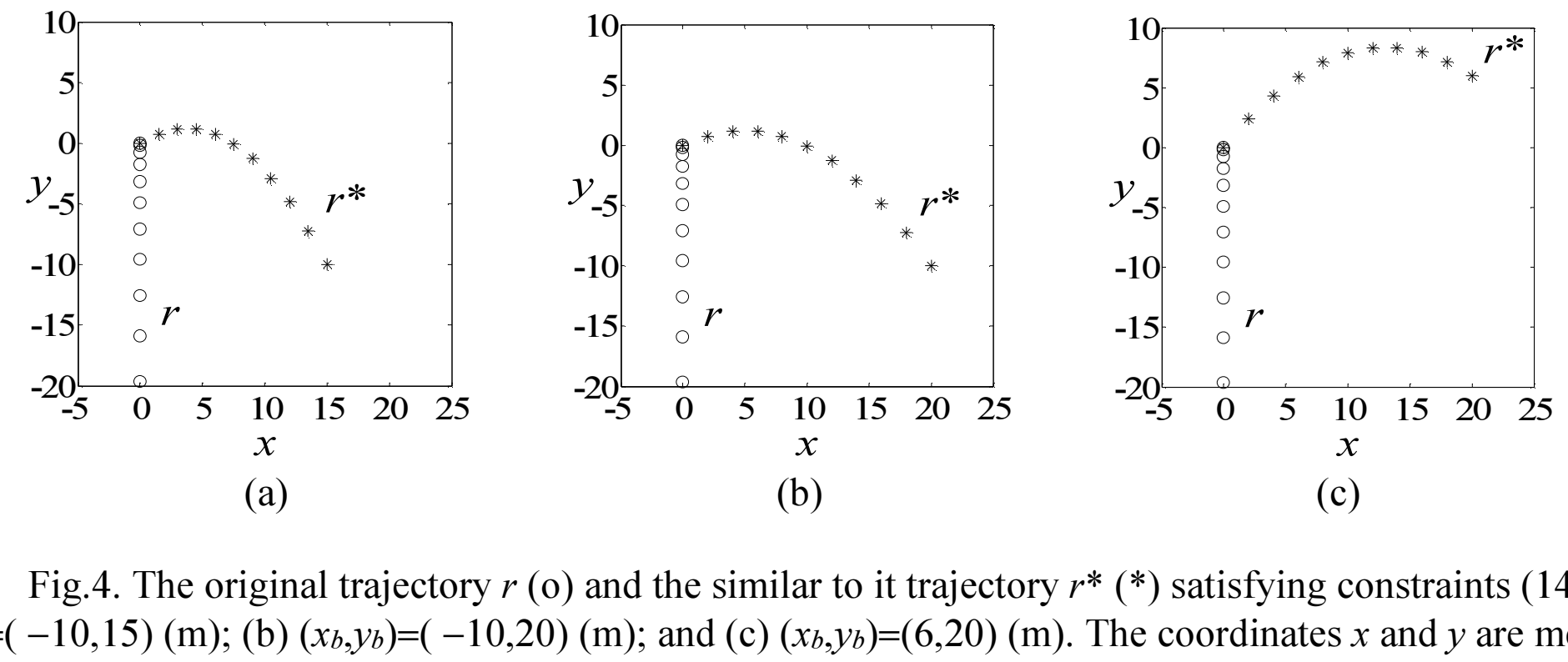

Fig.4. The original trajectory $r$ (o) and the similar to it trajectory $r^*$ (*) satisfying constraints (14) for (a) $(x_b,y_b)=(-10,15)$ (m); (b) $(x_b,y_b)=(-10,20)$ (m); and (c) $(x_b,y_b)=(6,20)$ (m). The coordinates $x$ and $y$ are measured in meters (m).

**Appendix**

In this appendix a MATLAB code for solving Example 3(c) is presented. The variables `A_`, `c_`, and `L_` are used for $\overline{A}$, $\overline{c}$ and $\overline{L}$, while `xs`, `ys`, and `us` are used for $x^*$, $y^*$ and $u^*$. The variable `l` is used for $\lambda$. To define the needed vectors and matrices, first the corresponding vectors and matrices composed of zeros and having the required size are defined.

```
function main

    N=101; M=6; a=0; b=2*3.141593; h=(b-a)/(N-1);

    t=zeros(N,1); x=zeros(N,1); y=zeros(N,1); u=zeros(2*N,1);
    A=zeros(M,2*N); c=zeros(M,1); A_=zeros(2*N,2*N); c_=zeros(2*N,1); L_=zeros(2*N,2*N);
```

```
        for i=1:N
            t(i)=a+(i-1)*h;
            x(i)=sin(2*t(i)); u(i)=x(i);
            y(i)=(1-sin(t(i)))*sin(t(i)); u(N+i)=y(i);
        end

        Sx=0; Sy=0; Tx=0; Ty=0;
        for i=1:N
            Sx=Sx+x(i); Sy=Sy+y(i); Tx=Tx+t(i)*x(i); Ty=Ty+(t(i)-a)*(t(i)-b)*y(i);
        end

        for i=1:N
            A(1,i)=1; A(2,N+i)=1; A(3,i)=t(i); A(4,N+i)=(t(i)-a)*(t(i)-b);
        end
        A(5,1)=1; A(5,N)=-1; A(6,N+1)=1; A(6,2*N)=-1;

        c(1)=Sx+10; c(2)=Sy; c(3)=Tx; c(4)=Ty+100; c(5)=0; c(6)=0;

        for i=1:N
            A_(1,i)=A(1,i); A_(N+1,N+i)=A(2,N+i);
        end

        c_(1)=c(1); c_(N+1)=c(2);

        L_(1,1)=-1; L_(1,2)=1; L_(N,N-1)=1; L_(N,N)=-1;
        L_(N+1,N+1)=-1; L_(N+1,N+2)=1; L_(2*N,2*N-1)=1; L_(2*N,2*N)=-1;
        for i=2:N-1
            L_(i,i-1)=1; L_(i,i)=-2; L_(i,i+1)=1;
            L_(N+i,N+i-1)=1; L_(N+i,N+i)=-2; L_(N+i,N+i+1)=1;
        end

        H=inv(L_+A_); d=A_*u-c_;

        l=(A*H*A')\(A*u-c-A*H*d)*2;
        us=u-H*(A'*l/2+d);

        for i=1:N
            xs(i)=us(i); ys(i)=us(N+i);
        end

        plot(x,y,'o',xs,ys,'*');
end
```